\documentclass{article}

\usepackage{amsmath,amssymb}

\begin{document}

\title{Sigma Functions for Telescopic Curves
} 

\author{Takanori Ayano}

\date{}

\maketitle

\begin{abstract}
In this paper we consider a symplectic basis of the first cohomology 
group and the sigma functions for algebraic curves expressed by a canonical form
using a finite sequence $(a_1,\cdots,a_t)$ of positive integers whose
greatest common divisor is equal to one (Miura \cite{Miura}).  
The idea is to express a non-singular algebraic curve by affine equations of $t$ variables whose
orders at infinity are $(a_1,\cdots,a_t)$. 
We construct a symplectic basis of the first cohomology 
group and the sigma functions for telescopic curves, i.e., the curves such that the number of
defining equations is exactly $t-1$ in the Miura canonical form. 
The largest class of curves for which such construction has been obtained thus far is $(n,s)$-curves
(\cite{Buch3}\cite{Nakayashiki}), which are telescopic because they are expressed in the Miura canonical form with $t=2$, 
$a_1=n$, and $a_2=s$, and the number of defining equations is one.

\end{abstract}

\section{Introduction}

Recently the theory of Abelian functions is attracting increasing interest in mathematical physics and applied mathematics. 
In particular the sigma functions for algebraic curves have been studied actively. 
In this paper we construct sigma functions explicitly for a class of algebraic curves for which such construction has not been obtained thus far. 

Let $C$ be a compact Riemann surface of genus $g$ and $H^1(C,\mathbb{C})$ the first cohomology group, which is defined by the linear space of second kind differentials modulo meromorphic exact forms. 
We say a meromorphic differential on $C$ to be second kind if it is locally exact. 

We consider a basis of $H^1(C,\mathbb{C})$ consisting of $\dim_{\mathbb{C}} H^1(C,\mathbb{C})=2g$ elements (cf. \cite{Lang}, pp.29-31, Theorem 8.1,8.2) . 
In particular, in order to construct sigma functions explicitly, we wish to construct a basis (symplectic basis) $\{du_i,dr_i\}_{i=1}^g$ of $H^1(C,\mathbb{C})$ such that 

\begin{enumerate}

\item $du_i$ is holomorphic on $C$ for each $i$, and 

\item $du_i\circ du_j=dr_i\circ dr_j=0$ and $du_i\circ dr_j=\delta_{ij}$ for each $i,j$, 

\end{enumerate}

\noindent where the operator $\circ$ is the intersection form on $H^1(C,\mathbb{C})$ defined by 
\[\eta\circ\eta'=\sum_p\mbox{Res}(\int^p\eta)\eta'(p)\]
for second kind differentials $\eta,\eta'$ (the summation is over all the singular points of $\eta$ and $\eta'$, and Res means taking a residue at a point).  

In order to express defining equations of $C$, we use a canonical form for expressing non-singular algebraic curves introduced by Miura \cite{Miura}. 
Given a finite sequence $(a_1,\dots,a_t)$ of positive integers whose
greatest common divisor is equal to one, Miura \cite{Miura} introduced a non-singular algebraic curve determined by the sequence $(a_1,\dots,a_t)$. 
The idea is to express a non-singular algebraic curve by affine equations of $t$ variables whose
orders at infinity are $(a_1,\cdots,a_t)$. 
Any non-singular algebraic curve is birationally equivalent to a curve expressed in the Miura canonical form (cf. \cite{Miura}).

Klein \cite{Klein1}\cite{Klein2} extended the elliptic sigma functions to the case of hyperelliptic curves of genus $g$, which are expressed in the Miura canonical form with $t=2$, $a_1=2$, and $a_2=2g+1$. 
Bukhshtaber et al. \cite{Buch3} and Nakayashiki \cite{Nakayashiki} extended Klein's sigma functions to the case of more general plane algebraic curves called $(n,s)$-curves, which are expressed in the Miura canonical form with $t=2$, $a_1=n$, and $a_2=s$. 
In this paper we give an explicit construction of sigma functions for telescopic curves, i.e., the curves such that the number of defining equations is exactly $t-1$ in
the Miura canonical form.
The telescopic curves contain the $(n,s)$-curves as special cases. 
Recently Matsutani \cite{Matsutani} constructed sigma functions for $(3,4,5)$-curves, which are not telescopic. 

The plan of this paper is as follows. 
In section 2 we recall the definition of the Miura canonical form.
In section 3 we construct the holomorphic 1-forms $\{du_i\}_{i=1}^g$ for the telescopic curves. 
In section 4 we construct the second kind differentials $\{dr_i\}_{i=1}^g$ for the telescopic curves and show that the set $\{du_i,dr_i\}_{i=1}^g$ is a symplectic basis of the first cohomology group. 
In section 5 we construct sigma functions for the telescopic curves. 

Throughout this paper, $\mathbb{N},\mathbb{N}_{+},\mathbb{Z}$, and $\mathbb{C}$ denote the set of non-negative integers, positive integers, integers, and complex numbers, respectively.

\section{Miura canonical form}

Miura [13] introduced a canonical form of defining equations for any non-singular algebraic curve. 
Here we recall the definition of the Miura canonial form. 

Let $t\ge2$, $a_1,\dots,a_t$ positive integers such that GCD$\{a_1,\dots,a_t\}=1$, 
$A_t=(a_1,\dots,a_t)\in\mathbb{N}_{+}^t$, and $\langle A_t\rangle=a_1\mathbb{N}+\cdots+a_t\mathbb{N}$, assuming that the order of $a_1,\dots,a_t$ is fixed. 
For the map $\Psi:\mathbb{N}^t \to \langle A_t\rangle$ defined by $\Psi((m_1,\dots,m_t))=\sum_{i=1}^ta_im_i$, 
we define the order $<$ in $\mathbb{N}^t$ so that $M<M'$ for $M=(m_1,\dots,m_t)$ and $M'=(m_1',\dots,m_t')$ if 

\begin{enumerate}
\item $\Psi(M)<\Psi(M')$ or 

\item $\Psi(M)=\Psi(M')$ and $m_1=m_1',\dots, m_{i-1}=m_{i-1}', m_i>m_i'$ for some $i \;(1\le i\le t)$. 

\end{enumerate}

\noindent Let $M(a)$ be the minimum element with respect to the order $<$ in $\mathbb{N}^t$ satisfying $\Psi(M)=a\in\langle A_t\rangle$. 
We define $B(A_t)\subseteq \mathbb{N}^t$ and $V(A_t)\subseteq\mathbb{N}^t\backslash B(A_t)$ by 
\[B(A_t)=\{M(a)\;|\;a\in \langle A_t\rangle\}\]
\noindent and  
\[V(A_t)\!=\!\{L\in {{\mathbb N}}^t\backslash B(A_t)\;|\;L=M+N,M\in {{\mathbb N}}^t\backslash B(A_t), N\in {\mathbb N}^t\Rightarrow N=(0,\dots,0)\},\]
respectively. 

Hereafter $\mathbb{C}[X]:=\mathbb{C}[X_1,\dots,X_t]$ denotes the polynomial ring over $\mathbb{C}$ of $t$-variables $X_1,\dots,X_t$. 
For $A\subset\mathbb{C}[X]$, Span$\{A\}$ and $(A)$ denote the linear space over $\mathbb{C}$ generated by $A$ and the ideal in $\mathbb{C}[X]$ generated by $A$, respectively. 
Also $X^M, M=(m_1,\dots,m_t)$, denotes $X^M=X_1^{m_1}\cdots X_t^{m_t}$ for simplicity. 

For $M\in V(A_t)$ we define the polynomial $F_M(X)\in\mathbb{C}[X]$ by 
\begin{equation}
F_M(X)=X^M-X^L-\sum_{\{N\in B(A_t)|\Psi(N)<\Psi(M)\}}\lambda_NX^N,\;\;\;\;\;\;\lambda_N\in\mathbb{C},\label{last}
\end{equation}
where $L$ is the element of $B(A_t)$ satisfying $\Psi(L)=\Psi(M)$. 
We assume that the set of polynomials $\{F_M\;|\;M\in V(A_t)\}$ satisfies the following condition:
\begin{equation}
\mbox{Span}\{X^N\;|\;N\in B(A_t)\}\cap (\{F_M\;|\;M\in V(A_t)\})=\{0\}.\label{cd}
\end{equation}

Let $I= (\{F_M\;|\;M\in V(A_t)\})$, $R=\mathbb{C}[X]/I$, $x_i$ the image of $X_i$ for the projection $\mathbb{C}[X]\to R$, and $K$ the total quotient ring of $R$.  
Then we have the following three propositions. 
Because there exists no paper where proofs are written in English, we give complete proofs in Appendix.

\vspace{2ex}

{\bf Proposition 2.1} (Miura \cite{Miura}). 

\vspace{1ex}

{\itshape \noindent (i) The set $\{x^N\;|\;N\in B(A_t)\}$ is a basis of $R$ over $\mathbb{C}$, where $x=(x_1,\dots,x_t)$.

\vspace{1ex}

\noindent (ii) The ring $R$ is an integral domain, therefore $K$ is the quotient field of $R$. 
 
\vspace{1ex}

\noindent (iii) The field $K$ is an algebraic function field of one variable over $\mathbb{C}$. 

\vspace{1ex}

\noindent (iv) There exists a discrete valuation $v_{\infty}$ of $K$ such that $(x_i)_{\infty}=a_iv_{\infty}$ for any $i$, 
where $(x_i)_{\infty}$ denotes the pole divisor of $x_i$ (cf.\cite{Stichtenoth} p.19). 

}

\vspace{2ex}

Let $C^{\scriptsize\mbox{aff}}=\{(z_1,\dots,z_t)\in \mathbb{C}^t\;|\;f(z_1,\dots,z_t)=0, \; \forall f\in I\}$. 
From Proposition 2.1 (ii) (iii), $C^{\scriptsize\mbox{aff}}$ is an affine algebraic curve in $\mathbb{C}^t$.  
Hereafter we assume that $C^{\scriptsize\mbox{aff}}$ is non-singular. 
For $k\in \mathbb{N}$ we define $L(kv_{\infty})=\{f\in K\;|\;(f)+kv_{\infty}\ge0\}\cup\{0\}$, where $(f)$ denotes the divisor of $f$, i.e., $(f)=\sum_vv(f)\cdot v$.   

\vspace{2ex}

{\bf Proposition 2.2} (Miura \cite{Miura}). 
{\itshape (i) $R=\bigcup_{k=0}^{\infty}L(kv_{\infty})$. 

\vspace{1ex}

\noindent (ii) The map $\phi$
\[C^{\scriptsize\mbox{aff}}\to \{\mbox{discrete valuation of}\;\;K\}\backslash \{v_{\infty}\}\]
\[p\to v_p\hspace{28ex}\]
is bijective, where $v_p$ is the discrete valuation corresponding to $p\in C^{\scriptsize\mbox{aff}}$ (cf. \cite{Silverman}, p.21,22).}

\vspace{2ex}

Let $C$ be the compact Riemann surface corresponding to $C^{\scriptsize\mbox{aff}}$. 
From Proposition 2.2 (ii), $C$ is obtained from $C^{\scriptsize\mbox{aff}}$ by adding one point, say $\infty$, where the discrete valuation corresponding to $\infty$ is $v_{\infty}$. 
It is known that any non-singular algebraic curve is birationally equivalent to such $C$ for some $A_t$ (cf. \cite{Miura}). 
Hereafter we represent each curve $C$ by the sequence $A_t=(a_1,\dots,a_t)$ and call $(a_1,\dots,a_t)$-curve for short. 
 
The sequence $A_t=(a_1,\dots,a_t)$ is called telescopic if for any $i$ ($2\le i\le t$) 
\[\frac{a_i}{d_i}\in \frac{a_1}{d_{i-1}}\mathbb{N}+\cdots+\frac{a_{i-1}}{d_{i-1}}\mathbb{N},\;\;\;\;d_i:=\mbox{GCD}\{a_1,\cdots,a_i\}.\]
Note that $A_2=(a_1,a_2)$ is always telescopic. 

\vspace{2ex}

{\bf Proposition 2.3} (Miura \cite{Miura}). 
{\itshape If $A_t$ is telescopic,  then the condition (\ref{cd}) is satisfied and we have the following properties. 

\vspace{1ex}

\noindent (i) $B(A_t)=\{(m_1,\dots,m_t)\in\mathbb{N}^t\;|\;0\le m_i\le d_{i-1}/d_i-1,\;2\le i\le t\}$.

\vspace{1ex}

\noindent (ii) $V(A_t)=\{(d_{i-1}/d_i)\;{\mathbf e}_i\;|\;2\le i\le t\}$, where ${\mathbf e}_i$ is the $i$-th unit vector in ${\mathbb Z}^t$.

\vspace{1ex}

\noindent (iii) The genus $g$ of $C$ is  
\begin{equation}
g=\frac{1}{2}\left\{(1-a_1)+\sum_{i=2}^t\left(\frac{d_{i-1}}{d_i}-1\right)a_i\right\}.\label{genus}
\end{equation}}

\vspace{2ex}

Note that $\sharp V(A_t)$ is the number of defining equations, where $\sharp$ denotes the number of elements. 
From Lemma C.1 (iv) in Appendix, we obtain $\sharp V(A_t)\ge t-1$. 
If $A_t$ is telescopic, then from Proposition 2.3 (ii) we obtain $\sharp V(A_t)=t-1$. 
On the other hand Suzuki \cite{Suzuki} proved that if $\sharp V(A_t)=t-1$, then $A_t$ is telescopic by rearranging the elements in a proper order.

From Proposition 2.3, the defining equations of a telescopic $(a_1,\dots,a_t)$-curve are given as follows: for $2\le i\le t$,  
\[F_i(X_1,\dots,X_t)=X_i^{d_{i-1}/d_i}-\prod_{j=1}^tX_j^{m_{ij}}-\sum\lambda_{j_1,\dots,j_t}^{(i)}X_1^{j_1}\cdots X_t^{j_t},\]
where $(m_{i1},\dots,m_{it})\in B(A_t)$ such that $\sum_{j=1}^t a_jm_{ij}=a_id_{i-1}/d_i$, $\lambda_{j_1,\dots,j_t}^{(i)}\in\mathbb{C}$, and the sum is over all $(j_1,\dots,j_t)\in B(A_t)$ such that $\sum_{k=1}^ta_kj_k<a_id_{i-1}/d_i$.    
Assign degrees as 
\[\deg X_k=a_k,\;\;\;\;\deg \lambda_{j_1,\dots,j_t}^{(i)}=a_id_{i-1}/d_i-\sum_{k=1}^ta_kj_k.\]

\vspace{2ex}

EXAMPLE 1. $A_2=(n,s)$, $\;\;n,s\in\mathbb{N}_{+}$, $\;\;\mbox{GCD}\{n,s\}=1$.

\vspace{1ex}

\noindent Since $A_2=(n,s)$ is telescopic, from Proposition 2.3 (ii), we have $V(A_2)=\{(0,n)\}$. 
Therefore we have 
\[F_2(X_1,X_2)=X_2^n-X_1^s-\sum_{nj_1+sj_2<ns}\lambda^{(2)}_{j_1,j_2}X_1^{j_1}X_2^{j_2},\]
which is the $(n,s)$-curve introduced in \cite{Buch2}. 
In particular we obtain the elliptic curves if $n=2$ and $s=3$ and the hyperelliptic curves of genus $g$ if $n=2$ and $s=2g+1$. 

\vspace{2ex}

EXAMPLE 2. $A_3=(4,6,5)$.

\vspace{1ex}

\noindent Since $A_3=(4,6,5)$ is telescopic, from Proposition 2.3 (ii), we have $V(A_3)=\{(0,2,0), (0,0,2)\}$. 
Therefore we have 
\[F_2(X_1,X_2,X_3)=X_2^2-X_1^3-\lambda^{(2)}_{0,1,1}X_2X_3-\lambda^{(2)}_{1,1,0}X_1X_2-\lambda^{(2)}_{1,0,1}X_1X_3-\lambda^{(2)}_{2,0,0}X_1^2\]
\[-\lambda^{(2)}_{0,1,0}X_2-\lambda^{(2)}_{0,0,1}X_3-\lambda^{(2)}_{1,0,0}X_1-\lambda^{(2)}_{0,0,0}\]
and
\[F_3(X_1,X_2,X_3)=X_3^2-X_1X_2-\lambda^{(3)}_{1,0,1}X_1X_3-\lambda^{(3)}_{2,0,0}X_1^2-\lambda^{(3)}_{0,1,0}X_2-\lambda^{(3)}_{0,0,1}X_3\]
\[-\lambda^{(3)}_{1,0,0}X_1-\lambda^{(3)}_{0,0,0}.\]

\section{Holomorphic 1-forms for telescopic curves }

\hspace{4ex}Let $C$ be a telescopic $(a_1,\dots,a_t)$-curve and $\Gamma(C,\Omega_C^1)$ the linear space consisting of holomorphic 1-forms on $C$. 
In this section we construct a basis of $\Gamma(C,\Omega_C^1)$. 
Let $G$ be the matrix defined by 
\[G:=\begin{pmatrix} 
\displaystyle{\frac{\partial F_2}{\partial X_1}} & \dots & \displaystyle{\frac{\partial F_2}{\partial X_t}} \\
\hdotsfor{3}\\
\displaystyle{\frac{\partial F_t}{\partial X_1}} & \dots & \displaystyle{\frac{\partial F_t}{\partial X_t}}  
\end{pmatrix}\]
and $G_i$ the matrix obtained by removing the $i$-th column from $G$. 
Then we have the following theorem. 

\vspace{2ex}

{\bf Theorem 3.1.}
{\itshape The set 
\[P:=\left\{\displaystyle{\frac{x_1^{k_1}\cdots x_t^{k_t}}{\det G_1(x)}dx_1\;|\;(k_1,\dots,k_t)\in B(A_t), \;0\le \sum_{i=1}^ta_ik_i\le 2g-2}\right\}\]
is a basis of $\Gamma(C,\Omega_C^1)$ over $\mathbb{C}$, 
where $\det G_1(x)$ denotes $\det G_1(X=x)$. 
 }

\vspace{2ex}

We order the elements of $P$ in the ascending order with respect to the order at $\infty$ and  write $\{du_1,\dots,du_g\}$. 

\vspace{2ex}

In order to prove Theorem 3.1, we need some lemmas. 

\vspace{2ex}

{\bf Lemma 3.1.}
{\itshape If  $\det G_i(p)\neq 0$ for $p=(p_1,\dots,p_t)\in C^{\scriptsize\mbox{aff}}$ and $1\le i\le t$, then $v_p(x_i-p_i)=1$. }

\vspace{2ex}

Proof. Without loss of generality, we assume $i=1$. 
Suppose $v_p(x_1-p_1)\ge 2$. 
Then there exists $k$ ($2\le k\le t$) such that $v_p(x_k-p_k)=1$. 
In fact, if $v_p(x_k-p_k)\ge2$ for any $k$, then $v_p(f)\ge2$ or $v_p(f)=0$ for any $f\in R$. 
Then $v_p(g)\ge2$ or $v_p(g)=0$ for any $g\in R_p$, where $R_p$ is the localization of $R$ at $p$.
This contradicts that $R_p$ is a discrete valuation ring.  

There exist $\{\gamma_{ij}, \delta_{j_1,\dots,j_t}^{(i)}\}\in\mathbb{C}$ such that for $2\le i\le t$ 
\[F_i(X_1,\dots,X_t)=\sum_{j=1}^t \gamma_{ij}(X_j-p_j)+\sum_{j_1+\dots +j_t\ge 2} \delta_{j_1,\dots,j_t}^{(i)}(X_1-p_1)^{j_1}\cdots(X_t-p_t)^{j_t},\]
where $\gamma_{ij}=\frac{\partial F_i}{\partial X_j}(p)$. 
Since $F_i(x_1,\dots,x_t)=0$ and $v_p(x_1-p_1)\ge 2$, we have $v_p\left(\sum_{j=2}^t \gamma_{ij}(x_j-p_j)\right)=v_p\left((x_k-p_k)(\sum_{j=2}^t\gamma_{ij}\frac{x_j-p_j}{x_k-p_k})\right)\ge2$.  
Since $v_p(x_k-p_k)=1$, we have $\sum_{j=2}^t\gamma_{ij}b_j=0$, where $b_j=\left(\frac{x_j-p_j}{x_k-p_k}\right)(p)$. 
Therefore we obtain 
\[G_1(p)\begin{pmatrix}
b_2\\
\cdot\\
\cdot\\
b_t
\end{pmatrix}=
\begin{pmatrix}
0\\
\cdot \\
\cdot \\
0
\end{pmatrix}.\]
Since $b_k=1(\neq0)$, we have $\det G_1(p)=0$. 
This contradicts the assumption of Lemma 3.1. 
Therefore we obtain $v_p(x_1-p_1)=1$. 

\hspace{75ex}$\square$

\vspace{2ex}

{\bf Lemma 3.2.}
\noindent {\itshape (i) As an element of $K$, we have $\det G_1(x)\neq0$. 

\vspace{1ex}

\noindent (ii) $\displaystyle{\mbox{div}\left(\frac{dx_1}{\det G_1(x)}\right)=(2g-2)\infty}.$}

\vspace{1ex}

Proof. 
Since the differential $d\left(F_i(x_1,\dots,x_t)\right)=0$ for any $i$, we have 
\[
G(x)\begin{pmatrix}
dx_1\\
\cdot \\
\cdot \\
dx_t 
\end{pmatrix}=
\begin{pmatrix}
0\\
\cdot \\
\cdot \\
0
\end{pmatrix}\;\;.\label{A}
\]
By multiplying some elementary matrices on the left, the above equation becomes
\[\begin{pmatrix}
w_2\; z_{22} \;z_{23} \cdot \cdot \cdot z_{2t}\\
w_3\;\; 0\;\;\; z_{33} \cdot \cdot \cdot z_{3t}\\
\dots\\
w_t\;\; 0\;\; \cdot \cdot \cdot \;\;\;\;\;\;z_{tt}
\end{pmatrix}
\begin{pmatrix}
dx_1\\
\cdot\\
\cdot\\
dx_t
\end{pmatrix}=
\begin{pmatrix}
0\\
\cdot \\
\cdot \\
0
\end{pmatrix}.\]
Since $C^{\scriptsize\mbox{aff}}$ is non-singular, for any $p\in C^{\scriptsize\mbox{aff}}$ there exists $i$ such that $\det G_i(p)\neq0.$ 
Therefore we have $w_t\neq0$ or $z_{tt}\neq0$ as elements of $K$. 
Since $v_{\infty}(x_j)=-a_j$, we have $x_j\notin \mathbb{C}$, therefore $dx_j\neq0$ for any $j$. 
Since $w_tdx_1=z_{tt}dx_t$, we have $w_t\neq0$ and $z_{tt}\neq0$. 
Therefore, by multiplying some elementary matrices on the left, the above equation becomes
\[\begin{pmatrix}
w_2'\; z_{22} \;z_{23} \cdot \cdot \cdot 0\\
w_3'\;\; 0\;\;\; z_{33} \cdot \cdot \cdot 0\\
\dots\\
w_t\;\; 0\;\;\;\;\cdot \cdot \;\;\;\;z_{tt}
\end{pmatrix}
\begin{pmatrix}
dx_1\\
\cdot\\
\cdot\\
dx_t
\end{pmatrix}=
\begin{pmatrix}
0\\
\cdot \\
\cdot \\
0
\end{pmatrix}.\]
Similarly we obtain
\[
\begin{pmatrix}
w_2''\;\; z_{22} \;\;0 \;\cdot \cdot \cdot\; 0\\
w_3''\;\; 0\;\;\; z_{33} \cdot \cdot \cdot \;0\\
\dots\\
w_t''\;\; 0\;\;\; \cdot \cdot \cdot \;\;\;\;z_{tt}
\end{pmatrix}
\begin{pmatrix}
dx_1\\
\cdot\\
\cdot\\
dx_t
\end{pmatrix}=
\begin{pmatrix}
0\\
\cdot \\
\cdot \\
0
\end{pmatrix},\label{D}
\]
where $w_2'',\dots,w_t'', z_{22},\dots,z_{tt}\in K$ are non-zero. 
Therefore we obtain $\det G_1(x)=\pm z_{22}\cdots z_{tt}\neq0$, which complete the proof of (i). 

Next we prove that the 1-form $dx_1/\det G_1(x)$ is both holomorphic and non-vanishing on $C^{\scriptsize\mbox{aff}}$. 
When $\det G_1(p)\neq0$ for $p\in C^{\scriptsize\mbox{aff}}$, from Lemma 3.1, $dx_1/\det G_1(x)$ is both holomorphic and non-vanishing at $p$. 
Suppose $\det G_1(p)=0$ for $p\in C^{\scriptsize\mbox{aff}}$. 
Since $C^{\scriptsize\mbox{aff}}$ is non-singular, there exists $i\;(2\le i\le t)$ such that $\det G_i(p)\neq 0$.  
Since $w_i''dx_1+z_{ii}dx_i=0$, we have $w_i''z_{22}\cdots \widehat{z_{ii}}\cdots z_{tt}dx_1+z_{22}\cdots z_{tt}dx_i=0$, where $\widehat{z_{ii}}$ denotes to remove $z_{ii}$.  
Therefore we obtain 
\[(-1)^{i-2}\det G_i(x)dx_1+\det G_1(x)dx_i=0.\]  
Since $\det G_1(x)\neq0$ and $\det G_i(x)\neq0$,  we have 
\[\frac{dx_1}{\det G_1(x)}=(-1)^{i-1}\frac{dx_i}{\det G_i(x)}.\]
Therefore, from $\det G_i(p)\neq0$ and Lemma 3.1, $dx_1/\det G_1(x)$ is holomorphic and non-vanishing at $p$. 
On the other hand, by Riemann-Roch's theorem, we have $\deg\mbox{div}(dx_1/\det G_1(x))=2g-2$, which complete the proof of (ii). 

\hspace{75ex}$\square$

\vspace{2ex}

Proof of Theorem 3.1.
From Lemma 3.2 and Proposition 2.1 (i), we have $P\subset \Gamma(C,\Omega_C^1)$ and the elements of $P$ are linearly independent. 
Since $\dim_{\mathbb{C}}\Gamma(C,\Omega_C^1)=g$,  it is sufficient to prove $\sharp P=g$. 
It is well-known that there are $g$ gap values at $\infty$ from 0 to $2g-1$. 
Since $\dim_{\mathbb{C}}L((2g-1)v_{\infty})=\dim_{\mathbb{C}}L((2g-2)v_{\infty})=g$ (Riemann-Roch's theorem),  $2g-1$ is a gap value at $\infty$. 
Therefore, from Proposition 2.1 (i) and Proposition 2.2 (i), we have $\sharp\{(k_1,\dots,k_t)\in B(A_t) \;|\; 0\le \sum_{i=1}^ta_ik_i\le 2g-2\}=g$, which complete the proof of Theorem 3.1. 

\hspace{75ex}$\square$

\section{Second kind differentials for telescopic curves}

In this section we construct $dr_i$ for a telescopic $(a_1,\dots,a_t)$-curve $C$. 
For $2\le i\le t$ and $1\le j\le t$, 
let \[h_{ij}=\frac{F_i(Y_1,\dots,Y_{j-1},X_j,X_{j+1},\dots,X_t)-F_i(Y_1,\dots,Y_{j-1},Y_j,X_{j+1},\dots,X_t)}{X_j-Y_j}\]
and 
\[H=\begin{pmatrix} 
h_{22} & \dots & h_{2t} \\
\hdotsfor{3}\\
h_{t2} & \dots & h_{tt}  
\end{pmatrix}.\]
We consider the 1-form 
\[\Omega(x,y):=\frac{\det H(x,y)}{(x_1-y_1)\det G_1(x)}dx_1\]
and the bilinear form (cf. \cite{Nakayashiki}, p.181, 2.4)
\begin{equation}
\hat{\omega}(x,y):=d_y\Omega(x,y)+\sum c_{i_1,\dots,i_t;j_1,\dots,j_t}\frac{x_1^{i_1}\cdots x_t^{i_t}y_1^{j_1}\cdots y_t^{j_t}}{\det G_1(x)\det G_1(y)}dx_1dy_1\label{cccc}
\end{equation}
on $C\times C$, where $x=(x_1,\dots,x_t)$, $y=(y_1,\dots,y_t)$, $c_{i_1,\dots,i_t;j_1,\dots,j_t}\in\mathbb{C}$, $(i_1,\dots,i_t)\in B(A_t)$ satisfying $0\le \sum_{k=1}^ta_ki_k\le 2g-2$, and $(j_1,\dots,j_t)\in B(A_t)$. 

We take a basis $\{\alpha_i,\beta_i\}_{i=1}^g$ of the homology group $H_1(C,\mathbb{Z})$ such that their intersection numbers are $\alpha_i\circ \alpha_j=\beta_i\circ\beta_j=0$ and $\alpha_i\circ\beta_j=\delta_{ij}$. 

\vspace{2ex}

DEFINITION 4.1. (cf. \cite{Nakayashiki}, p.181, 2.4)
Let $\Delta=\{(p,p)\;|\;p\in C\}$. 
A meromorphic symmetric bilinear form $\omega(x,y)$ on $C\times C$ is called a normalized fundamental form if the following conditions are satisfied. 

\vspace{1ex}

\noindent (i) $\omega(x,y)$ is holomorphic except $\Delta$ where it has a double pole. 
For $p\in C$ take a local coordinate $s$ around $p$. 
Then the expansion in $s(x)$ at $s(y)$ is of the form
\[\omega(x,y)=\left(\frac{1}{(s(x)-s(y))^2}+regular\right)ds(x)ds(y).\]

\noindent (ii) $\displaystyle{\int_{\alpha_i}\omega=0}$ for any $i$, where the integration is with respect to any one of the variables.

\vspace{2ex}

\noindent Normalized fundamental form exists and unique (cf. \cite{Nakayashiki} p.182). 
Then we have the following theorem. 

\vspace{2ex}

{\bf Theorem 4.1.}

{\itshape \noindent (i) There exists a set of $c_{i_1,\dots,i_t;j_1,\dots,j_t}$ such that $\hat{\omega}(x,y)=\hat{\omega}(y,x)$, non-zero $c_{i_1,\dots,i_t;j_1,\dots,j_t}$ is a homogeneous polynomial of $\{\lambda_{l_1,\dots,l_t}^{(i)}\}$ of 
degree 
\[2\sum_{k=2}^t\frac{d_{k-1}}{d_k}a_k-\sum_{k=1}^t(i_k+j_k+2)a_k,\] 
and $c_{i_1,\dots,i_t;j_1,\dots,j_t}=0$ if $2\sum_{k=2}^t\frac{d_{k-1}}{d_k}a_k-\sum_{k=1}^t(i_k+j_k+2)a_k<0$. 

\vspace{1ex}

\noindent For a set of $c_{i_1,\dots,i_t;j_1,\dots,j_t}$ such that $\hat{\omega}(x,y)=\hat{\omega}(y,x)$, we have the following properties. 

\vspace{1ex}

\noindent (ii) The bilinear form $\hat{\omega}$ satisfies the condition (i) of Definition 4.1. 

\vspace{1ex}

\noindent (iii) For $du_i:=(x_1^{k_{i1}}\cdots x_t^{k_{it}}/\det G_1(x))dx_1$, we define 
\[dr_i=\sum_{j_1,\dots,j_t}c_{k_{i1},\dots,k_{it};j_1,\dots,j_t}\frac{y_1^{j_1}\cdots y_t^{j_t}}{\det G_1(y)}dy_1.\]
Then $dr_i$ is a second kind differential for any $i$, and the set $\{du_i,dr_i\}_{i=1}^g$ is a symplectic basis of $H^1(C,\mathbb{C})$. 

 }

\vspace{2ex}


Let $B$ be the set of branch points for the map $x_1:C\to\mathbb{P}^1, (x_1,\dots,x_t)\to[x_1:1]$ (cf. \cite{Silverman}, p.24, Example 2.2). 
Since the ramification index of the map $x_1$ at $\infty$ is $a_1$, we have $\deg x_1=a_1$
(cf. \cite{Silverman}, p.28, Proposition 2.6).  
For $p\in C$ we set $x_1^{-1}(x_1(p))=\{p^{(0)},p^{(1)},\dots,p^{(a_1-1)}\}$ with $p=p^{(0)}$, 
where the same $p^{(i)}$ is listed according to its ramification index. 

\vspace{2ex}

{\bf Lemma 4.1.}
{\itshape Let $U$ be a domain in $\mathbb{C}$, $f(z_1,z_2)$ a holomorphic function on $U\times U$, and $g(z)=f(z,z)$.
If $g\equiv0$ on $U$, then there exists a holomorphic function $h(z_1,z_2)$ on $U\times U$ such that $f(z_1,z_2)=(z_1-z_2)h(z_1,z_2)$. }

\vspace{1ex}

Proof.
Let $h(z_1,z_2)=f(z_1,z_2)/(z_1-z_2)$. 
Given $z_1$, $h(z_1,\cdot\;)$ has a singularity only at $z_1$, where its singularity is removable. 
Therefore $h(z_1,\cdot\;)$ is holomorphic on $U$. 
Similarly $h(\;\cdot,z_2)$ is holomorphic on $U$. 
Therefore $h$ is holomorphic on $U\times U$. 

\hspace{75ex}$\square$

\vspace{2ex}

{\bf Lemma 4.2.}
{\itshape The 1-form $\Omega(x,y)$ is holomorphic except $\Delta\cup\{(p^{(i)},p)\;|\;i\neq0, \;p\in B \;or\; p^{(i)}\in B\}\cup C\times\{\infty\}\cup\{\infty\}\times C$. 
}

\vspace{1ex}

Proof.
Since $dx_1/\det G_1(x)$ is holomorphic on $C$ (cf. Lemma 3.2), 
$\Omega(x,y)$ is holomorphic except $\Delta\cup\{(p^{(i)},p)\;|\;p\in C, i\neq0\}\cup C\times\{\infty\}\cup\{\infty\}\times C$. 
We prove that $\Omega(x,y)$ is holomorphic on $\{(p^{(i)},p)\;|\;i\neq0, p\notin B,  p^{(i)}\notin B\}$. 
We have 
\begin{equation}
F_i(X_1,\dots,X_t)=\sum_{j=1}^th_{ij}\cdot(X_j-Y_j)+F_i(Y_1,\dots,Y_t).\label{Fi}
\end{equation}
Set $X=x$ and $Y=y$, then we have 
\[\sum_{j=1}^th_{ij}(x,y)\cdot(x_j-y_j)=0.\]
Take $(p^{(i)},p)\in C\times C$ such that $i\neq0, p\notin B$, and $p^{(i)}\notin B$, then we have 
\[\begin{pmatrix} 
h_{21} & \dots & h_{2t} \\
\hdotsfor{3}\\
h_{t1} & \dots & h_{tt}  
\end{pmatrix}_{X=p^{(i)}, Y=p}
\begin{pmatrix} 
p_1^{(i)}-p_1 \\
\cdot \\
p_t^{(i)}-p_t 
\end{pmatrix}
=\begin{pmatrix} 
0 \\
\cdot \\
0
\end{pmatrix}.\]
Since $p_1^{(i)}-p_1=0$,  we have 
\[H(p^{(i)}, p)\begin{pmatrix}
p_2^{(i)}-p_2 \\
\cdot \\
p_t^{(i)}-p_t
\end{pmatrix}
=\begin{pmatrix}
0 \\ 
\cdot \\
0
\end{pmatrix}.\]
Since $(p_2^{(i)}-p_2,\dots,p_t^{(i)}-p_t)\neq(0,\dots,0)$, we have 
$\det H(p^{(i)}, p)=0$. 
Since $p\notin B$ and $p^{(i)}\notin B$, we can take $(x_1,y_1)$ as a local coordinate around $(p^{(i)},p)$.  
Therefore, from Lemma 4.1, there exists a holomorphic function $h(x_1,y_1)$ around $(p^{(i)}, p)$ such that $\det H(x,y)=(x_1-y_1)h(x_1,y_1)$. 
Therefore $\Omega(x,y)$ is holomorphic at $(p^{(i)},p)$. 

\hspace{75ex}$\square$

\vspace{2ex}

{\bf Lemma 4.3.}
{\itshape Let $p\notin B$, $s$ a local coordinate around $p$. 
Then the expansion of $\Omega(x,y)$ in $s(y)$ at $s(x)$ is of the form 
\[\Omega(x,y)=\left(\frac{-1}{s(y)-s(x)}+regular\right)ds(x).\]}

Proof.
Set $Y=y$ in (\ref{Fi}), then we have
\[F_i(X_1,\dots,X_t)=\sum_{j=1}^th_{ij}(X,y)\cdot(X_j-y_j).\]
Therefore we obtain 
\[\frac{\partial F_i}{\partial X_k}(x_1,\dots,x_t)=\sum_{j=1}^t\frac{\partial h_{ij}}{\partial X_k}(x,y)\cdot(x_j-y_j)+h_{ik}(x,y).\]
Set $x=y$, then we have 
\[\frac{\partial F_i}{\partial X_k}(x_1,\dots,x_t)=h_{ik}(x,x).\]
Therefore we obtain $\det G_1(x)=\det H(x,x)$. 
On the other hand, since $p\notin B$, we can take $(x_1,y_1)$ as a local coordinate around $(p,p)$. 
Since $p\notin B$, we have $\det G_1(p)\neq0$. 
In fact, if $\det G_1(p)=0$, then $dx_1/\det G_1(x)$ is not holomorphic at $p$, which contradicts Lemma 3.2 (ii). 
Therefore $\det H(x,y)/\det G_1(x)$ is holomorphic at $(p,p)$. 
Therefore, from Lemma 4.1, there exists a holomorphic function $\tilde{h}(x_1,y_1)$ around $(p,p)$ such that $\det H(x,y)/\det G_1(x)=1+(x_1-y_1)\tilde{h}(x_1,y_1)$. 
Therefore we obtain Lemma 4.3. 

\hspace{75ex}$\square$

\vspace{2ex}

{\bf Lemma 4.4.}
{\itshape When we express 
\[\det H(X,Y)=\sum\epsilon_{m_1,\dots,m_t,n_1,\dots,n_t}X_1^{m_1}\cdots X_t^{m_t}Y_1^{n_1}\cdots Y_t^{n_t},\]
we have $\sum_{k=1}^ta_k(m_k+n_k)\le\sum_{k=2}^ta_k\left((d_{k-1}/d_k)-1\right)$. }

\vspace{2ex}

Proof.
When we express 
\[F_i(X_1,\dots,X_t)=\sum_{k=0}^m\tilde{F}_{ik}^{(j)}(X_1,\dots,X_{j-1},X_{j+1},\dots,X_t)X_j^k,\]
we have
$h_{ij}=\sum_{k=1}^m\tilde{F}_{ik}^{(j)}(Y_1,\dots,Y_{j-1},X_{j+1},\dots,X_t)\sum_{l=0}^{k-1}X_j^lY_j^{k-l-1}$. 
Assign degrees as $\deg Y_k=a_k$, then $h_{ij}$ is a homogeneous polynomial of $\{\lambda_{j_1,\dots,j_t}^{(i)}, X_k, Y_k\}$ of degree $a_id_{i-1}/d_i-a_j$. 
Therefore we obtain Lemma 4.4. 

\hspace{75ex}$\square$

\vspace{2ex}

{\bf Lemma 4.5.}
{\itshape The meromorphic bilinear form $d_y\Omega(x,y)$ is holomorphic except $\Delta\cup\{(p^{(i)},p)\;|\;i\neq0,\; p\in B\;or\;p^{(i)}\in B\}\cup C\times\{\infty\}$. }

\vspace{1ex}

Proof.
It is sufficient to prove that $d_y\Omega(x,y)$ is holomorphic at $(\infty, y)$, $y\neq\infty$. 
From Lemma 4.4,  with respect to $x$, we obtain 
\[v_{\infty}\left(\det H(x,y)\right)\ge-\sum_{k=2}^ta_k\left((d_{k-1}/d_k)-1\right).\] 
If $v_{\infty}\left(\det H(x,y)\right)>-\sum_{k=2}^ta_k\left((d_{k-1}/d_k)-1\right)$, 
then from Lemma 3.2 (ii) and Proposition 2.3 (iii) we obtain $v_{\infty}\left(\Omega(x,y)\right)\ge0$. 
Therefore $d_y\Omega(x,y)$ is holomorphic at $(\infty, y)$. 
If $v_{\infty}\left(\det H (x,y)\right)=-\sum_{k=2}^ta_k\left((d_{k-1}/d_k)-1\right)$,
then $v_{\infty}\left(\Omega(x,y)\right)=-1$. 
Let $s$ be a local coordinate around $\infty$, then from Lemma 4.4 there exists a constant $e$ (which does not depend on $y$) such that
\[\Omega(x,y)=\left(\frac{e}{s}+\mbox{regular}\right)ds.\]
Therefore $d_y\Omega(x,y)$ is holomorphic at $(\infty,y)$, $y\neq\infty$. 

\hspace{75ex}$\square$

\vspace{2ex}

{\bf Lemma 4.6.}
{\itshape Let $\omega$ be the normalized fundamental form. 
Then there exist second kind defferentials $d\hat{r}_i\;(1\le i\le g)$ which are holomorphic except $\{\infty\}$
and satisfy the equation
\[\omega(x,y)-d_y\Omega(x,y)=\sum_{i=1}^gdu_i(x)d\hat{r}_i(y).\]}

Proof.
Set $B_2=\{(p^{(i)},p)\;|\;p\in B\backslash\{\infty\}\;\mbox{or}\;p^{(i)}\in B\backslash\{\infty\}\}$ in the proof of \cite{Nakayashiki} Lemma 5, then proof of Lemma 4.6 is similar to that of \cite{Nakayashiki} Lemma 5. 

\hspace{75ex}$\square$

\vspace{2ex}

{\bf Lemma 4.7.}
{\itshape Let $Q$ be the linear space consisting of meromorphic differentials on $C$ which are singular only at $\infty$ and 
\[S=\left\{(x_1^{i_1}\cdots x_t^{i_t}/\det G_1(x))dx_1\;|\;(i_1,\dots,i_t)\in B(A_t)\right\}.\] 
Then $S$ is a basis of $Q$. } 

\vspace{1ex}

Proof.
For $\eta\in Q$ 
we consider the meromorphic function $\eta/\frac{dx_1}{\det G_1(x)}$. 
From Lemma 3.2 (ii), it may have a pole only at $\infty$. 
From Proposition 2.1 (i) and Proposition 2.2 (i), $\eta/\frac{dx_1}{\det G_1(x)}$ is a linear combination of $x_1^{i_1}\cdots x_t^{i_t}$ with $(i_1,\dots,i_t)\in B(A_t)$ and the elements of $S$ are linearly independent. 
 
\hspace{75ex}$\square$

\vspace{2ex}

Proof of Theorem 4.1 (i).  We have 
\[d_y\Omega(x,y)=\frac{\{\sum_{k=1}^t(-1)^{k+1}(x_1-y_1)\frac{\partial \det H}{\partial Y_k}(x,y)\det G_k(y)\}+\det G_1(y)\det H(x,y)}{(x_1-y_1)^2\det G_1(x)\det G_1(y)}dx_1dy_1.\]
Then, $\det G_k$, $\det H$, and $(\partial\det H/\partial Y_k)$ are homogeneous polynomials of $\{\lambda_{j_1,\dots,j_t}^{(i)},X_j,Y_j\}$ of degree $\sum_{i=2}^t\frac{d_{i-1}}{d_i}a_i-\sum_{i\neq k}a_i$, 
$\sum_{i=2}^t(\frac{d_{i-1}}{d_i}-1)a_i$, and $\{\sum_{i=2}^t(\frac{d_{i-1}}{d_i}-1)a_i\}-a_k$, respectively. 
Let us write
\[d_y\Omega(x,y)=\frac{\sum q_{i_1,\dots,i_t;j_1,\dots,j_t}x_1^{i_1}\cdots x_t^{i_t}y_1^{j_1}\cdots y_t^{j_t}}{(x_1-y_1)^2\det G_1(x)\det G_1(y)}dx_1dy_1,\]
where $(i_1,\dots,i_t), (j_1,\dots,j_t)\in B(A_t)$, and $q_{i_1,\dots,i_t;j_1,\dots,j_t}\in\mathbb{C}$. 
Then $q_{i_1,\dots,i_t;j_1,\dots,j_t}\in\mathbb{Z}[\{\lambda_{l_1,\dots,l_t}^{(i)}\}]$ and $q_{i_1,\dots,i_t;j_1,\dots,j_t}$ is homogeneous of degree $2\sum_{k=2}^t(\frac{d_{k-1}}{d_k}-1)a_k-\sum_{k=1}^t(i_k+j_k)a_k$. 
Note that if $(m_1,\dots,m_t)\in B(A_t)$, then $(m_1+m,m_2,\dots,m_t)\in B(A_t)$ for $m\in\mathbb{N}$. 
Therefore we obtain 
\[\sum c_{i_1,\dots,i_t;j_1,\dots,j_t}\frac{x_1^{i_1}\cdots x_t^{i_t}y_1^{j_1}\cdots y_t^{j_t}}{\det G_1(x)\det G_1(y)}\]
\[=\frac{\sum(c_{i_1-2,\dots,i_t;j_1,\dots,j_t}-2c_{i_1-1,\dots,i_t;j_1-1,\dots,j_t}+c_{i_1,\dots,i_t;j_1-2,\dots,j_t})x_1^{i_1}\cdots x_t^{i_t}y_1^{j_1}\cdots y_t^{j_t}}{(x_1-y_1)^2\det G_1(x)\det G_1(y)},\]
where $(i_1,\dots,i_t), (j_1,\dots,j_t)\in B(A_t)$. 
Therefore $\hat{\omega}(x,y)=\hat{\omega}(y,x)$ is equivalent to 
\[c_{i_1-2,\dots,i_t;j_1,\dots,j_t}-2c_{i_1-1,\dots,i_t;j_1-1,\dots,j_t}+c_{i_1,\dots,i_t;j_1-2,\dots,j_t}-c_{j_1-2,\dots,j_t;i_1,\dots,i_t}\]
\[+2c_{j_1-1,\dots,j_t;i_1-1,\dots,i_t}-c_{j_1,\dots,j_t;i_1-2,\dots,i_t}=q_{j_1,\dots,j_t;i_1,\dots,i_t}-q_{i_1,\dots,i_t;j_1,\dots,j_t}.\]
By Lemma 4.6, 4.7, the system of the above linear equations has a solution. 
Moreover it has a solution such that each $c_{i_1,\dots,i_t;j_1,\dots,j_t}$ is a linear combination of $q_{i_1',\dots,i_t';j_1',\dots,j_t'}$ satisfying $i_1'+j_1'=i_1+j_1+2$, $(i_k',j_k')=(i_k,j_k)$ or $(i_k',j_k')=(j_k,i_k)$ 
for $k=2,\dots,t$. 
In particular one can take $c_{i_1,\dots,i_t;j_1,\dots,j_t}$ such that $c_{i_1,\dots,i_t;j_1,\dots,j_t}=0$ if $2\sum_{k=2}^t\frac{d_{k-1}}{d_k}a_k-\sum_{k=1}^t(i_k+j_k+2)a_k<0$ and 
\[\deg c_{i_1,\dots,i_t;j_1,\dots,j_t}=2\sum_{k=2}^t\frac{d_{k-1}}{d_k}a_k-\sum_{k=1}^t(i_k+j_k+2)a_k\]
if $c_{i_1,\dots,i_t;j_1,\dots,j_t}\neq0$.

\hspace{75ex}$\square$

\vspace{2ex}

Proof of Theorem 4.1 (ii).
From Lemma 4.6, $d_y\Omega(x,y)$ is holomorphic except $\Delta\cup C\times\{\infty\}$ and so is $\hat{\omega}$. 
Since $\hat{\omega}(x,y)=\hat{\omega}(y,x)$, $\hat{\omega}$ is holomorphic except $\Delta$. 
From the definition of $dr_i$, we obtain 
\[\hat{\omega}-\omega=\sum_{i=1}^gdu_i(x)(dr_i(y)-d\hat{r}_i(y)).\]
On the other hand $\hat{\omega}-\omega$ is holomorphic except $\Delta$ and $\sum_{i=1}^gdu_i(x)(dr_i(y)-d\hat{r}_i(y))$ is holomorphic except $C\times\{\infty\}$. 
Therefore $\hat{\omega}-\omega$ is holomorphic except $\{\infty\}\times\{\infty\}$.  
Therefore $\hat{\omega}-\omega$ and $dr_i-d\hat{r}_i$ are holomorphic on $C\times C$ and $C$ respectively, which complete the proof of Theorem 4.1 (ii). 

\hspace{75ex}$\square$

\vspace{2ex}

Proof of Theorem 4.1 (iii).
The 1-form $dr_i$ is a second kind differential. 
In fact $dr_i-d\hat{r}_i$ is holomorphic 1-form as is just proved in the proof of Theorem 4.1 (ii) and $d\hat{r}_i$ is a second kind differential from Lemma 4.6. 
Proof of Theorem 4.1 (iii) is similar to the case of the $(n,s)$-curves (cf. \cite{Nakayashiki} Lemma 7,8, Proposition 3). 

\hspace{75ex}$\square$

\section{Sigma functions for telescopic curves}

In this section we construct the sigma function for a telescopic $(a_1,\dots,a_t)$-curve $C$. 
First we take the following data. 

\begin{enumerate}
\item A basis $\{\alpha_i,\beta_i\}_{i=1}^g$ of the homology group $H_1(C,\mathbb{Z})$ such that their intersection numbers are $\alpha_i\circ \alpha_j=\beta_i\circ\beta_j=0$ and $\alpha_i\circ\beta_j=\delta_{ij}$. 

\item The symplectic basis $\{du_i,dr_i\}_{i=1}^g$ of the first cohomology group $H^1(C,\mathbb{C})$ constructed in section 3 and 4. 

\end{enumerate}

\noindent We define the period matrices by 
\[2\omega_1=\left(\int_{\alpha_j}du_i\right),\;2\omega_2=\left(\int_{\beta_j}du_i\right),\;-2\eta_1=\left(\int_{\alpha_j}dr_i\right),\;-2\eta_2=\left(\int_{\beta_j}dr_i\right).\]
Then $\omega_1$ is invertible. 
Set $\tau=\omega_1^{-1}\omega_2$,  
then $\tau$ is symmetric and Im $\tau>0$. 
By the Riemann's bilinear relation 
\[2\pi i\eta\circ\eta'=\sum_{i=1}^g\left(\int_{\alpha_i}\eta\int_{\beta_i}\eta'-\int_{\alpha_i}\eta'\int_{\beta_i}\eta\right),\]
the matrix 
\[M:=\begin{pmatrix}
\omega_1 & \omega_2 \\
\eta_1 & \eta_2
\end{pmatrix}
\]
satisfies 
\[M\begin{pmatrix}
0 & I_g \\
-I_g & 0
\end{pmatrix}
{}^tM=-\frac{\pi \sqrt{-1}}{2}
\begin{pmatrix}
0 & I_g \\
-I_g & 0
\end{pmatrix},\]
where $I_g$ denotes the unit matrix of degree $g$. 
Since $\eta_1\omega_1^{-1}$ is symmetric (cf. \cite{Nakayashiki} Lemma 8), we obtain the following proposition. 

\vspace{2ex}

{\bf Proposition 5.1.} (generalized Legendre relation) 
\[{}^tM\begin{pmatrix}
0 & I_g \\
-I_g & 0
\end{pmatrix}
M=-\frac{\pi \sqrt{-1}}{2}
\begin{pmatrix}
0 & I_g \\
-I_g & 0
\end{pmatrix}.\]

\hspace{75ex}$\square$

\vspace{2ex}

Let $\delta=\tau\delta'+\delta''$ be the Riemann's constant of $C$ with respect to our choice $(\infty, \{\alpha_i,\beta_i\}_{i=1}^g)$. 
Since the divisor of the holomorphic 1-form $du_g$ is $(2g-2)\infty$, the Riemann's constant $\delta$ becomes a half period. 
Then the sigma funtion $\sigma(u)$ associated with $C$ is defined as follows. 

\vspace{2ex}

DEFINITION 5.1.  (Sigma function)
For $u\in\mathbb{C}^g$ 
\[\sigma(u)=\sigma(u;M)=c\cdot\exp\left(\frac{1}{2}\;{}^tu\eta_1\omega_1^{-1}u\right)\theta\begin{bmatrix} \delta' \\ \delta'' \end{bmatrix}((2\omega_1)^{-1}u,\tau)\]
\[=c\cdot\exp\left(\frac{1}{2}\;{}^tu\eta_1\omega_1^{-1}u\right)\hspace{2ex}\]
\[\times\sum_{n\in\mathbb{Z}^g}\exp\left\{\pi \sqrt{-1}\;{}^t(n+\delta')\tau(n+\delta')+2\pi \sqrt{-1}\;{}^t(n+\delta')((2\omega_1)^{-1}u+\delta'')\right\},\]
where $c$ is a constant. 

\vspace{2ex}

\noindent By Proposition 5.1 we obtain the following proposition. 

\vspace{2ex}

{\bf Proposition 5.2.}
For any $m_1,m_2\in\mathbb{Z}^g$ and $u\in\mathbb{C}^g$, we have 
\[\sigma(u+2\omega_1m_1+2\omega_2m_2)/\sigma(u)=\exp\left(\pi \sqrt{-1}\;({}^tm_1m_2+2\;{}^t\delta'm_1-2\;{}^t\delta''m_2)\right)\]
\[\times\exp\left(\;{}^t(2\eta_1m_1+2\eta_2m_2)(u+\omega_1m_1+\omega_2m_2)\right).\]

\hspace{75ex}$\square$

\vspace{4ex}

REMARK. 
In this paper we have constructed sigma functions explicitly for telescopic curves. 
One can show that the first term of the series expansion around
the origin of the sigma functions for telescopic curves becomes Schur function corresponding
to the partition determined from the gap sequence at infinity and the expansion
coefficients are homogeneous polynomials of the coefficients of the defining equations
of the curve in a manner similar to \cite{Nakayashiki}.

\vspace{4ex}

\noindent {\Large \bf Appendix}

\appendix

\section{Proof of Proposition 2.1}

\hspace{4ex}{\bf Lemma A.1}. \hspace{3ex}{\itshape $V(A_t)+\mathbb{N}^t=\mathbb{N}^t\backslash B(A_t)$. }

\vspace{1ex}

Proof.
If $M\notin B(A_t)$ and $N\in\mathbb{N}^t$, then $M+N\notin B(A_t)$. 
Therefore we have $V(A_t)+\mathbb{N}^t\subset \mathbb{N}^t\backslash B(A_t)$.
Suppose $V(A_t)+\mathbb{N}^t\subsetneq\mathbb{N}^t\backslash B(A_t)$. 
Take $M_1\in\mathbb{N}^t\backslash B(A_t)$ satisfying $M_1\notin V(A_t)+\mathbb{N}^t$. 
Since $M_1\notin V(A_t)$ and $M_1\notin B(A_t)$, there exist $M_2\in\mathbb{N}^t\backslash B(A_t)$ and $(0,\dots,0)\neq N_1\in\mathbb{N}^t$ such that $M_1=M_2+N_1$. 
Since $M_1\notin V(A_t)+\mathbb{N}^t$, we have $M_2\notin V(A_t)+\mathbb{N}^t$.
Similarly, for the element $M_i\in\mathbb{N}^t\backslash B(A_t)$ satisfying $M_i\notin V(A_t)+\mathbb{N}^t$, there exist $M_{i+1}$ and $N_i$ such that $M_{i+1}\in\mathbb{N}^t\backslash B(A_t)$, $M_{i+1}\notin V(A_t)+\mathbb{N}^t$, $(0,\dots,0)\neq N_i\in\mathbb{N}^t$, and $M_i=M_{i+1}+N_i$. 
Therefore there exists a infinite sequence $\Psi(M_1)>\Psi(M_2)>\cdots>\Psi(M_i)>\cdots$.
This is contradiction. 

\hspace{75ex}$\square$

\vspace{2ex}

Proof of Proposition 2.1 (i).
From (\ref{cd}) it is sufficient to prove 
\[\mbox{Span}\{X^N\;|\;N\in B(A_t)\}+ (\{F_M\;|\;M\in V(A_t)\})=\mathbb{C}[X].\]
We prove that for any $T\in\mathbb{N}^t$
\[X^T\in\mbox{Span}\{X^N\;|\;N\in B(A_t)\}+ (\{F_M\;|\;M\in V(A_t)\})\]
by transfinite induction with respect to the well-order $<$ in $\mathbb{N}^t$. 
The statement is correct for the minimal element $T=(0,\dots,0)$.  
Suppose that it is correct for any $U\in\mathbb{N}^t$ satisfying $U<T$. 
Since it is correct for $T\in B(A_t)$,  
we assume $T\notin B(A_t)$. From Lemma A.1, there exist $M\in V(A_t)$ and $Z\in\mathbb{N}^t$ such that $T=M+Z$. 
Then we have $X^T=X^MX^Z=(X^M-F_M)X^Z+F_MX^Z$. 
For any monomial $X^U$ in $(X^M-F_M)X^Z$, we have $U<T$. 
Therefore, by the assumption of transfinite induction, the statement is correct for $T\notin B(A_t)$. 

\hspace{75ex}$\square$

\vspace{2ex}

\noindent We define the function $o : R\to \mathbb{N}\cup\{-\infty\}$ by 

\[
o(f)=\left\{
\begin{array}{l}
-\infty\;\;\;\;\;\;\;\;\;\;\;\;\;\;\;\;\;\;\;\;\;\;\;\;\;\;\;\;\;\mbox{for}\;\;f=0 \\ 
\max\{\Psi(N)\;|\;\lambda_N\neq0\}\;\;\mbox{for}\;\;f\neq0
\end{array}
\right.,
\]
where for $f\neq0$ we express $f=\sum_N\lambda_Nx^N$ with $\lambda_N\in\mathbb{C}$ and $N\in B(A_t)$. 

\vspace{2ex}

{\bf Lemma A.2}. 
{\itshape $o(x^T)=\Psi(T)$ for any $T\in\mathbb{N}^t$. }

\vspace{1ex}

Proof.
We prove the statement by transfinite induction with respect to the well-order $<$ in $\mathbb{N}^t$.
It is correct for the minimal element $T=(0,\dots,0)\in\mathbb{N}^t$. 
Suppose that it is correct for any $U\in\mathbb{N}^t$ satisfying $U<T$.
Since it is correct for $T\in B(A_t)$,  we assume $T\notin B(A_t)$. 
From Lemma A.1, there exist $M\in V(A_t)$ and $Z\in\mathbb{N}^t$ such that $T=M+Z$. 
Then we have $X^T=X^MX^Z=(X^M-F_M)X^Z+F_MX^Z$.
Since $X^M-F_M=X^L+\sum_N\lambda_NX^N$ from (\ref{last}),  we have $x^T=(x^L+\sum_N\lambda_Nx^N)x^Z=x^{L+Z}+\sum_N\lambda_Nx^{N+Z}$. 
Since $N+Z<L+Z<T$, 
by the assumption of transfinite induction, 
we have $o(x^{L+Z})=\Psi(L+Z)$ and $o(x^{N+Z})=\Psi(N+Z)$. 
Since $o(f+g)=\max\{o(f),o(g)\}$ for $f, g\in R$ satisfying $o(f)\neq o(g)$, 
we have $o(x^T)=o(x^{L+Z}+\sum_N\lambda_Nx^{N+Z})=o(x^{L+Z})=\Psi(L+Z)=\Psi(T). $

\hspace{75ex}$\square$

\vspace{2ex}
 
{\bf Lemma A.3.}
{\itshape The function $o$ satisfies the following properties: 

\vspace{1ex}

\noindent (i)\;\;$o(f)=-\infty$ if and only if $f=0$, 

\vspace{1ex}

\noindent (ii)\;\;$o(fg)=o(f)+o(g)$ for any $f,g\in R$, 
where we define $-\infty+(-\infty)=a+(-\infty)=(-\infty)+a=-\infty$ for $a\in\mathbb{N}$, 

\vspace{1ex}

\noindent (iii)\;\;$o(f+g)\le\max\{o(f),o(g)\}$, 

\vspace{1ex}

\noindent (iv)\;\;$o(R\backslash\{0\})=\langle A_t\rangle$, in particular $\mathbb{N}\backslash o(R\backslash\{0\})$ is a finite set,  and 

\vspace{1ex}

\noindent (v)\;\;$o(a)=0$ for any $0\neq a\in\mathbb{C}$. }

\vspace{2ex}

Proof. (i), (iii), (v), and $o(R\backslash\{0\})=\langle A_t\rangle$ are trivial.  
Since GCD$\{a_1,\dots,a_t\}=1$, $\mathbb{N}\backslash\langle A_t\rangle$ is a finite set (cf. \cite{inter}, Theorem 5). 
We prove (ii). 
If $f=0$ or $g=0$, then $o(fg)=o(f)+o(g)=-\infty$. 
Suppose $f\neq0$ and $g\neq0$. 
Then we can express 
\[f=\lambda_Mx^M+\sum_T\lambda_Tx^T\;\;\mbox{and}\;\;g=\tilde{\lambda}_Nx^N+\sum_Z\tilde{\lambda}_Zx^Z,\]
where $\lambda_M, \lambda_T, \tilde{\lambda}_N,\tilde{\lambda}_Z\in\mathbb{C}$, $\lambda_M\neq0,\tilde{\lambda}_N\neq0$, $M,T,N,Z\in B(A_t)$, $\Psi(T)<\Psi(M)$, and $\Psi(Z)<\Psi(N).$ 
From Lemma A.2, we have $o(fg)=o(\lambda_M\tilde{\lambda}_Nx^{M+N})=\Psi(M+N)=\Psi(M)+\Psi(N)=o(f)+o(g).$

\hspace{75ex}$\square$

\vspace{2ex}

Proof of Proposition (ii). 
Take $f,g\in R$ satisfying $fg=0$. 
Then, since $-\infty=o(fg)=o(f)+o(g)$, we have $o(f)=-\infty$ or $o(g)=-\infty$. 
Therefore we obtain $f=0$ or $g=0$. 

\hspace{75ex}$\square$

\vspace{2ex}

{\bf Lemma A.4}. 
{\itshape Let $B\subset\mathbb{N}^t$ be a set such that the restriction map of $\Psi : \mathbb{N}^t\to\langle A_t\rangle$ on $B$ is bijective. 
Then the set $\{x^M\;|\;M\in B\}\subset R$ is a basis of $R$ over $\mathbb{C}$. }

\vspace{1ex}

Proof.
Since $o(x^T)=\Psi(T)$ for $T\in\mathbb{N}^t$ and $o(f+g)=\max\{o(f), o(g)\}$ for $f, g\in R$ satisfying $o(f)\neq o(g)$, 
the elements of the set $\{x^M\;|\;M\in B\}$ are linearly independent. 
Since $R=\mbox{Span}\{x^N\;|\;N\in B(A_t)\}$, in order to prove $R=\mbox{Span}\{x^M\;|\;M\in B\}$, it is sufficient to prove $\mbox{Span}\{x^N\;|\;N\in B(A_t)\}\subset \mbox{Span}\{x^M\;|\;M\in B\}$. 
We prove  $\mbox{Span}\{x^N\;|\;N\in B(A_t),\;\Psi(N)\le m\}\subset \mbox{Span}\{x^M\;|\;M\in B,\;\Psi(M)\le m\}$ for any $m\in\mathbb{N}$ by induction. 
For $m=0$ the statement is trivial. 
Suppose that the statement is correct for any $i$ with $0\le i\le m-1$. 
If  $m\notin\langle A_t\rangle$, then since $\mbox{Span}\{x^M\;|\;M\in B,\;\Psi(M)\le m\}=\mbox{Span}\{x^M\;|\;M\in B,\;\Psi(M)\le m-1\}$ and $\mbox{Span}\{x^N\;|\;N\in B(A_t),\;\Psi(N)\le m\}=\mbox{Span}\{x^N\;|\;N\in B(A_t),\;\Psi(N)\le m-1\}$, 
the statement is correct. 
Suppose $m\in\langle A_t\rangle$. 
Take $T\in B$ satisfying $\Psi(T)=m$. 
If $T\in B(A_t)$,  then since $\mbox{Span}\{x^M\;|\;M\in B,\;\Psi(M)\le m\}=\mbox{Span}\{x^M\;|\;M\in B,\;\Psi(M)\le m-1\}\cup\mathbb{C}\{x^T\}$ and $\mbox{Span}\{x^N\;|\;N\in B(A_t),\;\Psi(N)\le m\}=\mbox{Span}\{x^N\;|\;N\in B(A_t),\;\Psi(N)\le m-1\}\cup\mathbb{C}\{x^T\}$, the statement is correct. 
Suppose $T\notin B(A_t)$. 
Then we can express $x^T=\lambda_Lx^L+\sum_N\lambda_Nx^N$, where 
$0\neq\lambda_L,\lambda_N\in\mathbb{C}$, $L,N\in B(A_t)$, $\Psi(L)=m$, and $\Psi(N)\le m-1$. 
Since 
$x^L=\lambda_L^{-1}(x^T-\sum_N\lambda_Nx^N)\in\mbox{Span}\{x^N\;|\;N\in B(A_t),\;\Psi(N)\le m-1\}\cup\mathbb{C}\{x^T\}\subset\mbox{Span}\{x^M\;|\;M\in B,\;\Psi(M)\le m-1\}\cup\mathbb{C}\{x^T\}\subset\mbox{Span}\{x^M\;|\;M\in B,\;\Psi(M)\le m\}$, 
we have $\mbox{Span}\{x^N\;|\;N\in B(A_t),\;\Psi(N)\le m\}\subset \mbox{Span}\{x^M\;|\;M\in B,\;\Psi(M)\le m\}$. 

\hspace{75ex}$\square$

\vspace{2ex}

{\bf Lemma A.5}. 
{\itshape Given $i$, there exists a set $T_i\subset\mathbb{N}^{i-1}\times\{0\}\times\mathbb{N}^{t-i}$ such that $\sharp T_i=a_i$ and for the set $B_i:=T_i+\{0\}^{i-1}\times\mathbb{N}\times\{0\}^{t-i}$ the restriction map of $\Psi : \mathbb{N}^t\to\langle A_t\rangle$ on $B_i$ is bijective. }

\vspace{1ex}

Proof.
Since GCD$\{a_1,\dots,a_t\}=1$, the set $\{c\in a_1\mathbb{N}+\cdots+a_{i-1}\mathbb{N}+a_{i+1}\mathbb{N}+\cdots+a_t\mathbb{N}\;|\;c\equiv\; j\; \mbox{mod}\;a_i\}$ is not empty for any $j$ with $0\le j\le a_i-1$. 
Let $c_j=\min\{c\in a_1\mathbb{N}+\cdots+a_{i-1}\mathbb{N}+a_{i+1}\mathbb{N}+\cdots+a_t\mathbb{N}\;|\;c\equiv\; j\; \mbox{mod}\;a_i\}$. 
Take $N_j\in\mathbb{N}^{i-1}\times\{0\}\times\mathbb{N}^{t-i}$ satisfying $\Psi(N_j)=c_j$. 
Let $T_i=\{N_j\;|\;0\le j\le a_i-1\}$. 
Then $T_i$ satisfies the conditions of Lemma A.5. 

\hspace{75ex}$\square$

\vspace{2ex}

Proof of Proposition 2.1 (iii).
Since $o(x^T)=\Psi(T)$ for $T\in\mathbb{N}^t$ and $o(f+g)=\max\{o(f), o(g)\}$ for $f, g\in R$ satisfying $o(f)\neq o(g)$,  the elements of the set $\{x^M\;|\;M\in\{0\}^{i-1}\times\mathbb{N}\times\{0\}^{t-i}\}\subset\mathbb{C}[x_i]$ are linearly independent. 
Therefore the extension of field $\mathbb{C}(x_i)/\mathbb{C}$ is a simple transcendental extension for any $i$.  
Next we prove $[K : \mathbb{C}(x_i)]\le a_i$ for any $i$. 
From Lemma A.4 and Lemma A.5, we have $R=\mathbb{C}[x_1,\dots,x_t]=\mbox{Span}\{x^M\;|\;M\in T_i+\{0\}^{i-1}\times\mathbb{N}\times\{0\}^{t-i}\}$. 
Therefore $\mathbb{C}[x_1,\dots,x_t]=\mathbb{C}[x_i]f_0+\cdots+\mathbb{C}[x_i]f_{a_i-1}$, where $f_j=x^{N_j}$ (see the proof of Lemma A.5 for $N_j$). 
Since $f_0=1$, we obtain the finite extension of integral domain $\mathbb{C}(x_i)\subset\mathbb{C}(x_i)f_0+\cdots+\mathbb{C}(x_i)f_{a_i-1}$.
Since $\mathbb{C}(x_i)$ is a field, $\mathbb{C}(x_i)f_0+\cdots+\mathbb{C}(x_i)f_{a_i-1}$ is also a field. 
Therefore we obtain $\mathbb{C}(x_i)f_0+\cdots+\mathbb{C}(x_i)f_{a_i-1}=K$ and $[K : \mathbb{C}(x_i)]\le a_i$. 

\hspace{75ex}$\square$

\vspace{2ex}

Proof of Proposition 2.1 (iv). 
We define the function $v_{\infty} : K\to\mathbb{Z}\cup\{\infty\}$ by 
\[
v_{\infty}(f)=\left\{
\begin{array}{l}
\infty\;\;\;\;\;\;\;\;\;\;\;\;\;\;\;\;\;\;\;\;\;\;\mbox{for}\;\;f=0 \\ 
-o(f_1)+o(f_2)\;\;\;\mbox{for}\;\;f\neq0
\end{array}
\right.,
\]
where for $f\neq0$ we express $f=f_1/f_2$ with $f_1,f_2\in R$. 
The definition of $v_{\infty}$ is well-defined. 
In fact, if $0\neq f=f_1/f_2=g_1/g_2$, then since $f_1g_2=g_1f_2\in R$, we have $o(f_1)+o(g_2)=o(f_1g_2)=o(g_1f_2)=o(g_1)+o(f_2)$. 
From Lemma A.3, one can check that the function $v_{\infty}$ is a discrete valuation of $K$. 
From Lemma A.2, we obtain $v_{\infty}(x_i)=-a_i$. 
From \cite{Stichtenoth} p.19 Theorem 1.4.11, we obtain $[K : \mathbb{C}(x_i)]=\deg (x_i)_{\infty}\ge \deg (a_iv_{\infty})=a_i$. 
On the other hand, in the proof of Proposition 2.1 (iii), we proved $[K : \mathbb{C}(x_i)]\le a_i$. 
Therefore we obtain $(x_i)_{\infty}=a_iv_{\infty}$. 

\hspace{75ex}$\square$

\section{Proof of Proposition 2.2}

\hspace{4ex}Proof of Proposition 2.2 (i).
It is trivial that $R\subset \bigcup_{k=0}^{\infty}L(kv_{\infty})$. 
On the other hand  we have 
\[\bigcup_{k=0}^{\infty}L(kv_{\infty})\subset\bigcap_{v\neq v_{\infty}}\mathcal{O}_v\subset\bigcap_{p\in C^{\scriptsize\mbox{aff}}}\mathcal{O}_p=R,\]
where $\mathcal{O}_v=\{f\in K\;|\;v(f)\ge0\}$ and $\mathcal{O}_p=\{f\in K\;|\;v_p(f)\ge0\}$ (see Proposition 2.2 (ii) for $v_p$). 

\hspace{75ex}$\square$

\vspace{2ex}

Proof of Proposition 2.2 (ii).
It is trivial that the map $\phi$ is injective. 
We prove that the map $\phi$ is surjective. 
Let $v$ be a discrete valuation such that $v\neq v_{\infty}$. 
Since $v(x_i)\ge0$ for any $i$, we have $R\subset \mathcal{O}_v$. 
Let $P$ be the maximal ideal of $\mathcal{O}_v$ and $m:=P\cap R$. 
Then we have 
\[\mathbb{C}\hookrightarrow R/m \hookrightarrow \mathcal{O}_v/P.\]
Since $[\mathcal{O}_v/P : \mathbb{C}]=1$, we have $\mathbb{C}\simeq R/m\simeq \mathcal{O}_v/P.$
Therefore $m$ is a maximal ideal. 
Let $R_m$ be the localization of $R$ with respect to $m$. 
Then $R_m$ and $\mathcal{O}_v$ are discrete valuation rings satisfying $R_m\subset\mathcal{O}_v$ and $P\cap R_m=mR_m. $
Therefore, from \cite{Hartshorne} p.40 Theorem 6.1A, we obtain $R_m=\mathcal{O}_v$. 
Since there exists $p\in C^{\scriptsize\mbox{aff}}$ such that $\mathcal{O}_p=R_m$, we have $\mathcal{O}_p=\mathcal{O}_v$. 
Therefore we obtain $v_p=v$ and the map $\phi$ is surjective. 

\hspace{75ex}$\square$

\section{Proof of Proposition 2.3}

\hspace{4ex}Let $T(A_t)=B(A_t)\cap\left(\{0\}\times\mathbb{N}^{t-1}\right).$ 

\vspace{2ex}

{\bf Lemma C.1}. 
{\itshape (i) $T(A_t)=\{M(b_i)\in B(A_t)\;|\;i=0,\dots,a_1-1\}$, where $b_i=\min\{b\in a_2\mathbb{N}+\cdots +a_t\mathbb{N}\;|\;b\equiv i\;\mbox{mod}\;a_1\}$. 
In particular $\sharp T(A_t)=a_1$. 

\vspace{1ex}

\noindent (ii) $B(A_t)=T(A_t)+\mathbb{N}\times\{0\}^{t-1}$. 

\vspace{1ex}

\noindent (iii) $V(A_t)\subset \{T(A_t)+{\mathbf e}_i\;|\;i=2,\dots,t\}\backslash T(A_t)\subset\{0\}\times\mathbb{N}^{t-1}$.

\vspace{1ex}

\noindent (iv) The set $\{0\}^{i-1}\times\mathbb{N}\times\{0\}^{t-i}\cap V(A_t)$ consists of only one element for any $i\;(2\le i\le t)$. }
 
\vspace{2ex}

Proof.
We have $M(b_i)=(m_1,\dots,m_t)\in \{0\}\times\mathbb{N}^{t-1}$. 
In fact,  if $m_1\neq0$, then we have 
$\Psi((0,m_2,\dots,m_t))\equiv b_i\equiv i\;\mbox{mod}\;a_1$ and $\Psi((0,m_2,\dots,m_t))<b_i$, 
which contradicts the definition of $b_i$. 
Therefore we have $M(b_i)\in T(A_t)$. 
For $M,N\in\{0\}\times\mathbb{N}^{t-1}$ satisfying $\Psi(M)>\Psi(N)$ and $\Psi(M)-\Psi(N)=ea_1$ for some $e\in\mathbb{N}_{+}$, we have $M\notin T(A_t)$. 
In fact, for $N':=(e,0,\dots,0)+N$, we have $M>N'$ and $\Psi(M)=\Psi(N')$, which means $M\notin B(A_t)$. 
Therefore we obtain (i). 

Next we prove $B(A_t)\subset T(A_t)+\mathbb{N}\times\{0\}^{t-1}$. 
Let $M=(m_1,\dots,m_t)\in B(A_t)$, $M_1=(0,m_2,\dots,m_t)$, and $M_2=(m_1,0,\dots,0)$. 
Since $M_1+M_2\in B(A_t)$, we have $M_1,M_2\in B(A_t)$. 
Since $M_1\in B(A_t)\cap(\{0\}\times\mathbb{N}^{t-1})=T(A_t)$, we have $M\in T(A_t)+\mathbb{N}\times\{0\}^{t-1}$.  
Suppose $B(A_t)\subsetneq T(A_t)+\mathbb{N}\times\{0\}^{t-1}.$ 
Then from (i) there exist $i \;(0\le i\le a_1-1)$ and $M_3\in\mathbb{N}\times\{0\}^{t-1}$ such that $M(b_i)+M_3\notin B(A_t)$. 
Take $N\in B(A_t)$ satisfying $\Psi(M(b_i)+M_3)=\Psi(N)$. 
Since $N\in B(A_t)\subset T(A_t)+\mathbb{N}\times\{0\}^{t-1}$ and $\Psi(N)\equiv i\;\mbox{mod}\;a_1$, there exists $M_4\in\mathbb{N}\times\{0\}^{t-1}$ such that $N=M(b_i)+M_4$. 
Therefore $M_3>M_4,\;M_3,M_4\in\mathbb{N}\times\{0\}^{t-1}$, and $\Psi(M_3)=\Psi(M_4)$, which is contradiction. 
Therefore we obtain $B(A_t)=T(A_t)+\mathbb{N}\times\{0\}^{t-1}.$

Next we prove $V(A_t)\subset\{0\}\times\mathbb{N}^{t-1}$.  
Let $M=(m_1,\dots,m_t)\in V(A_t),\;M_1=(0,m_2,\dots,m_t)$, and $M_2=(m_1,0,\dots,0)$. 
Since $M\notin B(A_t)$ and $M_2\in B(A_t)$, we have $M_1\notin B(A_t)$. 
From the definition of $V(A_t)$, we obtain $M_2=(0,\dots,0)$. 
Therefore we obtain $V(A_t)\subset\{0\}\times\mathbb{N}^{t-1}$. 

Let $M\in V(A_t)\subset\{0\}\times\mathbb{N}^{t-1}$. 
Since $M\neq(0,\dots,0)$, there exist $i\;(2\le i\le t)$ and $M_1\in\{0\}\times\mathbb{N}^{t-1}$ such that $M=M_1+{\mathbf e}_i$. 
Since $M_1\in B(A_t)$ from the definition of $V(A_t)$, 
we have $M_1\in B(A_t)\cap\left(\{0\}\times\mathbb{N}^{t-1}\right)=T(A_t)$. 
Therefore we obtain (iii). 

For $2\le i\le t$, the set $\{0\}^{i-1}\times\mathbb{N}\times\{0\}^{t-i}\cap\{\mathbb{N}^t\backslash B(A_t)\}$ is not empty. 
In fact, since 
\[\Psi((0,\dots,0,a_1,0,\dots,0))=\Psi((a_i,0,\dots,0))=a_1a_i,\]
we have $(0,\dots,0,a_1,0,\dots,0)>(a_i,0,\dots,0)$. 
Let $N_i$ be the minimal element of $\{0\}^{i-1}\times\mathbb{N}\times\{0\}^{t-i}\cap\{\mathbb{N}^t\backslash B(A_t)\}$. 
Then we obtain $\{0\}^{i-1}\times\mathbb{N}\times\{0\}^{t-i}\cap V(A_t)=\{N_i\}$. 
Therefore we obtain (iv). 

\hspace{75ex}$\square$

\vspace{2ex}

\noindent Let $SV(A_t)=\{N_i\;|\;2\le i\le t\}$ (see the proof of Lemma C.1 (iv) for $N_i$). 
For $F=\sum\lambda_NX^N\in\mathbb{C}[X]$, we define {\itshape multideg} of $F$ by
\[
\mbox{multideg}(F)=\left\{
\begin{array}{l}
-\infty\;\;\;\;\;\;\;\;\;\;\;\;\;\;\;\;\;\;\;\;\;\;\;\;\;\;\;\;\;\;\;\mbox{for}\;\;F=0 \\ 
\displaystyle{\max_{<}}\{N\in\mathbb{N}^t\;|\;\lambda_N\neq0\}\;\;\mbox{for}\;\;F\neq0
\end{array}
\right..
\]
Also we define {\itshape leading term} of $F$ by
\[
\mbox{LT}(F)=\left\{
\begin{array}{l}
0\;\;\;\;\;\;\;\;\;\;\;\;\;\mbox{for}\;\;F=0 \\ 
\lambda_TX^T\;\;\;\;\;\mbox{for}\;\;F\neq0,\;\;\mbox{where} \;\;T=\mbox{multideg}(F)
\end{array}
\right..
\]
For a ideal $J\subset \mathbb{C}[X]$, we define  
\[\Delta(J)=\mathbb{N}^t\backslash\bigcup_{F\in J\backslash\{0\}}\{\mbox{multideg}(F)+\mathbb{N}^t\}.\]
Then we have 
\begin{equation}
\mbox{Span}\{X^M\;|\;M\in\Delta(J)\}\cap J=\{0\}. \label{span}
\end{equation}

\vspace{1ex}

{\bf Lemma C.2}. 
{\itshape (i) $\{F_M\;|\;M\in SV(A_t)\}$ is a Gr\"{o}bner basis of the ideal $J:=\left(\{F_M\;|\;M\in SV(A_t)\}\right)$ with respect to the order $<$ in $\mathbb{N}^t$, 
i.e., $(\{\mbox{LT}(F)\;|\;F\in J\})=(\{\mbox{LT}(F_M)\;|\;M\in SV(A_t)\})$. 

\vspace{1ex}

\noindent (ii) Span$\{X^N\;|\;N\in B(A_t)\}\cap (\{F_M\;|\;M\in SV(A_t)\})=\{0\}$. }

\vspace{1ex}

Proof.
For $M,N\in SV(A_t)$ ($M\neq N$), we have L.C.M.$\{\mbox{LT}(F_M),\mbox{LT}(F_N)\}=\mbox{LT}(F_M)\mbox{LT}(F_N)$. 
Therefore, from \cite{Cox} p.102 Theorem 3 and p.103 Proposition 4, we obtain (i). 
From (i) we obtain $\Delta\left((\{F_M\;|\;M\in SV(A_t)\})\right)=\mathbb{N}^t\backslash\{SV(A_t)+\mathbb{N}^t\}\supset\mathbb{N}^t\backslash\{V(A_t)+\mathbb{N}^t\}=B(A_t)$, where the last equality is due to Lemma A.1. 
Since Span$\{X^N\;|\;N\in\Delta(\{F_M\;|\;M\in SV(A_t)\})\}\cap (\{F_M\;|\;M\in SV(A_t)\})=\{0\}$ from (\ref{span}), 
we have Span$\{X^N\;|\;N\in B(A_t)\}\cap (\{F_M\;|\;M\in SV(A_t)\})=\{0\}$.

\hspace{75ex}$\square$

\vspace{2ex}

{\bf Lemma C.3}. 
{\itshape If $A_t$ is telescopic, then the following properties are satisfied. 

\vspace{1ex}

\noindent (i) $T(A_t)=\{(0,m_2,\dots,m_t)\in\mathbb{N}^t\;|\;0\le m_i\le d_{i-1}/d_i-1, i=2,\dots,t\}$. 

\vspace{1ex}

\noindent (ii) $SV(A_t)=V(A_t)=\{(d_{i-1}/d_i)\;{\mathbf e}_i\;|\;2\le i\le t\}.$}

\vspace{1ex}

Proof.
Let $U=\{(0,m_2,\dots,m_t)\in\mathbb{N}^t\;|\;0\le m_i\le d_{i-1}/d_i-1, i=2,\dots,t\}$. 
Take $u=(0,u_2,\dots,u_t)\in U$ and $v=(0,v_2,\dots,v_t)\in U$ satisfying $u\neq v$. 
First we prove $\Psi(u)\not\equiv\Psi(v)$ mod $a_1$. 
Suppose that there exists an integer $w$ such that $\Psi(u)-\Psi(v)=wa_1$. 
Let $\rho$ be the positive integer such that $u_{\rho}\neq v_{\rho},\;u_{\rho+1}=v_{\rho+1},\dots,u_t=v_t$. 
Without loss of generality we assume $u_{\rho}>v_{\rho}$. 
Then we have $(u_{\rho}-v_{\rho})a_{\rho}=wa_1-\sum_{k=2}^{\rho-1}(u_k-v_k)a_k$ and $0<u_{\rho}-v_{\rho}< d_{\rho-1}/d_{\rho}$, which is contradiction. 
Therefore  we obtain $\Psi(u)\not\equiv\Psi(v)$ mod $a_1$. 
Since $A_t$ is telescopic, 
for any $u=(0,u_2,\dots,u_t)\in\mathbb{N}^t$, there exists $u'\in U$ such that $\Psi(u)\equiv\Psi(u')$ mod $a_1$. 
Since $\Psi(u)\ge\Psi(u')$ and $\sharp U=a_1$, we have $\{\Psi(u)\;|\;u\in U\}=\{b_0,\dots,b_{a_1-1}\}$, where $b_i=\min\{b\in a_2\mathbb{N}+\cdots +a_t\mathbb{N}\;|\;b\equiv i\;\mbox{mod}\;a_1\}$. 
Finally we prove $u\in B(A_t)$ for any $u\in U$. 
Take $u\in U$, then 
there exists $u''=(u''_1,\dots,u''_t)\in B(A_t)$ such that $\Psi(u)=\Psi(u'')$. 
Since $A_t$ is telescopic, we have $0\le u_j''<d_{j-1}/d_j$ for $2\le j\le t$. 
Since $u_1''=0$ from the definition of $b_i$, we obtain $u''\in U$. 
Therefore we obtain $u=u''\in B(A_t)$. 
From Lemma C.1 (i), we obtain (i). 
From Lemma C.1 (iii) (iv) and the definition of $V(A_t)$, we obtain (ii). 

\hspace{75ex}$\square$

\vspace{2ex}

Proof of Proposition 2.3.
From Lemma C.2 (ii) and Lemma C.3 (ii), the condition (\ref{cd}) is satisfied. 
From Lemma C.1 (ii) and Lemma C.3 (i), we obtain Proposition 2.3 (i). 
From Lemma C.3 (ii), we obtain Proposition 2.3 (ii). 
From Proposition 2.1 (i) and Proposition 2.2 (i), the gap values at $\infty$ are $\mathbb{N}\backslash\langle A_t\rangle$. 
Therefore, from \cite{inter} Theorem 5, we obtain Proposition 2.3 (iii). 

\hspace{75ex}$\square$

\vspace{4ex}

ACKNOWLEDGEMENT. 
The author would like to thank his supervisor Prof. Joe Suzuki for suggesting for the author extending the sigma functions by using the Miura canonical form. 
The author would like to thank Prof. Yoshihiro \^{O}nishi for his warm encouragements and valuable discussions. 
The author would like to thank Prof. Ryuichi Harasawa for valuable comments for the Miura canonical form. 
This research was supported by Grant-in-Aid for JSPS Fellows (22-2421) from Japan Society for the Promotion of Science.

\vspace{4ex}

\begin{flushright}

Takanori Ayano

Department of Mathematics 

Graduate School of Science

Osaka University

Machikaneyama-chou 1--1

Toyonaka-si,Osaka

560--0043, Japan

Email address: t-ayano@cr.math.sci.osaka-u.ac.jp

\end{flushright}

\end{document}